\documentclass{amsart}

\usepackage{graphicx}
\usepackage{pst-all}      
\usepackage{url}
\usepackage{float}
\usepackage{fancyvrb}

\begin{document}

\title{Another Diminimal Map on the Torus}


\author{Thom Sulanke}
\address{Department of Physics, Indiana University, Bloomington, 
	Indiana 47405 USA}
\email{tsulanke@indiana.edu}

\date{\today}

\begin{abstract}

This note adds one diminimal map on the torus to the published set of 55.  
It also raises to 15 the number of vertices for which all diminimal maps on 
the torus are known.

\end{abstract}

\maketitle

\section{Introduction}
\label{intro}

Riskin \cite{MR1164599}\cite{MR1317680} found two diminimal maps on the torus 
(denoted R1 and R2 in this note).
Henry \cite{henry1}\cite{henry}, using a computer found 53 more 
(H1 through H53).
All diminimal maps on the torus with fewer than 10 
vertices or fewer than 10 faces were found using this technique.
But Henry was not able to 
complete the generation process due to time constraints.

Using a completely different computer technique we generated all polyhedral 
maps on the torus with 15 or fewer vertices and checked each to see if it was 
diminimal.
While this technique found all diminimal maps on the torus with 15 or fewer 
vertices, there may be other diminimal maps on the torus with more than 
15 vertices.

\section{Definitions}
\label{definitions}

Following \cite{henry} we use these definitions.
We consider only graphs with no multiple edges and no loops.
If $G$ is a graph embedded in a surface, 
then the closure of each connected component of the complement of the graph is
called a {\em face}.
An embedded graph together with its embedding is called a {\em map} 
if each vertex 
of the graph has degree at least 3 and each face is a closed 2-cell.

If the intersection of two distinct faces of a map is empty, a vertex, 
or an edge then the faces are said to {\em meet properly}. 
If each face is simply connected and each 
pair of faces meet properly, then the map is a {\em polyhedral map}.
We refer to a polyhedral map on
the torus as a {\em toroidal polyhedral map} or TPM.

The operation {\em edge removing} is the process of obtaining one map from
another by removing a single edge. If the removal creates a vertex of degree 2
then the two
edges adjacent to that vertex are coalesced into one edge in the new map. 
Let $G$ be a
TPM and let $G'$ be the map obtained from $G$ by removing 
edge $e$ from $G$. 
If $G'$ is also a TPM, then edge $e$ is called 
{\em removable}.

The operation {\em edge shrinking} is the process of obtaining one map from
another by contracting an edge so that edge's two vertices become one vertex. 
If the contraction creates a two-sided face in the new map then the two edges 
of the two-sided face are replaced with a single edge. 
Let $G$ be a TPM and let $G'$ be the map obtained
from $G$ by shrinking edge $e$ in $G$. 
If $G'$ is also a TPM, then edge $e$ is 
called {\em shrinkable}.

If $G$ is a TPM with no shrinkable edges and no removable 
edges then $G$ is called {\em diminimal}.

\section{Generating TPMs}
\label{gentpm}

The author has developed the program {\em surftri} \cite{surftri} for 
generating embeddings on surfaces.
This program is based on the program {\em plantri} \cite{plantri} written by 
Brinkmann and McKay which generates planar graphs.
It is easy to restrict {\em surftri} to generate only {\em polyhedral maps}.
As a test of {\em surftri} we checked that all the known
diminimal TPMs were among the TPMs
with 15 or fewer vertices. 
The number of TPMs checked is shown in Table~\ref{Counts}.
About three CPU years were required to generate the TPMs 
with up to 15 vertices.
In the process we discovered an additional diminimal TPM 
(S1).
All the known diminimal TPMs are shown in Figs.~\ref{fig1}, \ref{fig2}, 
and~\ref{fig3}.

A list of the 56  known diminimal TPMs (H1, \ldots, H53, R1, R2, S1)
is in Figure~\ref{list}.
Each line of the list represents one diminimal TPM.
The format used is taken from {\em surftri} \cite{surftri} which generalizes 
the format in {\em plantri} \cite{plantri}.
The vertices are labeled by single letters.
The number is the number of vertices for the embedding.
For each vertex there is a list of its neighbors in cyclic order.
So 
$${\tt 7\ bcdefg,agdfec,abegfd,acfbge,adgcbf,aebdcg,afcedb}$$
 represents $K_7$ embedded on the torus.
This embedding has 7 vertices.  The vertex {\tt a} has 6 neighbors 
{\tt bcdefg} in cyclic order.
The vertex {\tt b} also has 6 neighbors {\tt agdfec} in cyclic order.
The final vertex {\tt g} has neighbors {\tt afcedb}.

\vfill

\begin{table}[b]
\centering
\begin{tabular}{r|r r}
 Vertices & TPMs & Diminimal TPMs \\
\hline
7    & 1 & 1 \\
8    & 33 & 2 \\
9    & 4713 & 11 \\
10   & 442429 & 19 \\
11   & 28635972 & 15 \\
12   & 1417423218 & 5 \\
13   & 58321972887 & 2 \\
14   & 2102831216406 & 1 \\
15   & 68781200467456 & 0 \\
\end{tabular}
\vspace{.1in}
\caption{Counts of TPMs and diminimal TPMs by number of vertices}
\label{Counts}
\end{table}

\begin{figure}
\begin{Verbatim}[fontsize=\footnotesize]
7 bcdefg,agdfec,abegfd,acfbge,adgcbf,aebdcg,afcedb
8 bcde,aefgdh,ahgf,afhbg,aghfb,behdcg,bfced,bdfec
8 bcd,aefgh,ahegf,afhge,bdgchf,behdcg,bfcedh,bgdfec
9 bcde,afghi,aifh,ahge,adihf,becg,bfid,bdcei,bhegc
9 bcd,aefgh,ahei,aihge,bdgcf,behi,biedh,bgdfc,cgfd
9 bcd,aefgh,ahig,afhe,bdhgi,bidg,bfceh,bgedic,chfe
9 bcd,aefgh,aheg,agfi,bicf,behdg,bfdcih,bgifc,dhge
9 bcd,aefgh,aig,afhe,bdhgif,beidg,bfceh,bgedi,chfe
9 bcde,afghi,aifd,achg,agihf,behcg,bfied,bdfei,bhegc
9 bcd,aefgh,ahegf,afhge,bdgci,bidcg,bfcedh,bgdic,ehf
9 bcd,aefgh,ahig,agfhe,bdhgi,bihdg,bfdceh,bgedfc,cfe
9 bcd,aefgh,ahfi,agfhe,bdhif,bechdg,bfdih,bgedfc,ceg
10 bcde,afgh,ahfi,aihg,agif,becj,bjed,bdjc,cejd,fhig
10 bcd,aefgh,ahij,ajhe,bdgi,bihj,bjeh,bgdfc,cfe,cgfd
10 bcd,aefg,ahei,aihj,bjcf,behi,bijh,cgdf,cjgfd,dgie
10 bcd,aefgh,ahij,afhe,bdgi,bidj,bjeh,bgdic,chfe,cgf
10 bcd,aefgh,ahig,afhe,bdji,bidg,bfcj,bjdic,chfe,ehg
10 bcd,aefgh,ahei,aihj,bjcf,behi,bijh,bgdfc,cgfd,dge
10 bcd,aefgh,aig,afhe,bdjif,beidg,bfcj,bjdi,chfe,ehg
10 bcd,aefgh,aig,afje,bdjif,behdg,bfcj,bjfi,che,dhge
10 bcd,aefg,ahi,ajge,bdgihj,bji,biedh,cgje,cegf,dfeh
10 bcd,aefgh,aig,afj,bjif,behdg,bfcjh,bgjfi,che,dhge
10 bcd,aefg,aghi,aije,bdgih,bhji,biejc,cjfe,cegfd,dfhg
10 bcd,aefgh,ahfi,agje,bdhif,becjg,bfdih,bgejc,ceg,dfh
10 bcd,aefgh,ahig,agfj,bjgi,bihdg,bfdceh,bgjfc,cfe,dhe
10 bcde,aefgh,ahif,ajhg,aijb,bjcg,bfid,bdjic,cheg,dfeh
10 bcd,aefg,aghi,afhj,bjih,bhdig,bfjhc,cgdfe,cejf,dgie
10 bcd,aefgh,ahfi,ajhe,bdhif,bechj,bjih,bgedfc,ceg,dgf
10 bcd,aefg,ahif,ajgi,bihj,bjcig,bfidh,cgje,cedgf,dfeh
10 bcd,aefgh,ahfi,agfj,bjif,bechdg,bfdih,bgjfc,ceg,dhe
10 bcd,aefgh,aij,afhge,bdgif,behdj,bjedh,bgdfi,che,cgf
11 bcd,aefg,ahf,aigj,bji,bick,bkdh,cgij,dfeh,dkhe,fjg
11 bcd,aefg,ahf,aij,bjhi,bick,bkih,cgej,dfeg,dkhe,fjg
11 bcd,aefgh,aig,afhe,bdjk,bkdg,bfcj,bjdi,chk,ehg,eif
11 bcd,aefgh,aij,afhe,bdgk,bkdj,bjeh,bgdi,chk,cgf,eif
11 bcd,aefgh,aig,afj,bjgk,bkdg,bfceh,bgji,chk,dhe,eif
11 bcd,aefgh,ahig,ajhe,bdki,bihj,bjck,bkdfc,cfe,dgf,ehg
11 bcd,aefg,aghi,aije,bdhk,bihj,bjkc,cejf,cfkd,dgfh,egi
11 bcd,aefgh,ahei,afhj,bjck,bkdi,bijh,bgdkc,cgf,dge,ehf
11 bcd,aefgh,aij,afke,bdgif,behdj,bjek,bkfi,che,cgf,dhg
11 bcd,aefg,ahi,ajk,bkif,behj,bjikh,cgkf,cegj,digf,dhge
11 bcd,aefg,ahi,ajk,bkhf,beji,bikjh,cgje,ckgf,dfhg,dgie
11 bcd,aefg,aghi,ajgk,bkh,bjig,bfkdhc,cgje,ckf,dfh,dgie
11 bcd,aefgh,aij,afhk,bkif,behdj,bjkh,bgdfi,che,cgf,dge
11 bcd,aefgh,aij,agfk,bkjf,beidg,bfdjh,bgki,chf,ceg,dhe
11 bcde,afgh,ahfi,aijk,akif,becj,bjik,bkjc,cegd,dgfh,dhge
12 bcd,aefg,ahf,aij,bki,bicl,blih,cgk,dfeg,dlk,ejh,fjg
12 bcd,aefg,ahi,ajkl,blj,bilk,bkjh,cgl,cfj,dieg,dgf,dfhe
12 bcd,aefg,ahi,ajk,bkhl,bikg,bfjh,cgek,cfl,dlg,dfhe,eji
12 bcd,aefg,ahi,ajk,bkhf,beji,bikl,cle,ckgf,dfl,dgie,gjh
12 bcd,aefgh,aij,akl,blif,bekj,bjlh,bgki,che,cgf,dfh,dge
13 bcd,aefg,ahi,ajk,blm,bikg,bfjh,cgl,cfm,dmg,dfl,ekh,eji
14 bcd,aef,agh,aij,bkl,bmn,cnk,cml,dln,dmk,ejg,eih,fhj,fig
9 bcde,afdg,ageh,ahbi,aicf,begi,bhfc,cigd,dfhe
9 bcde,afdg,agfh,ahbi,aihf,beci,bhic,cegd,dfge
13 bcde,afgh,ahij,ajkl,almf,beji,biml,blkc,ckgf,cfmd,dmih,dhge,egkj
\end{Verbatim}

\caption{List of Diminimal Maps on the Torus in {\em plantri} format}
\label{list}
\end{figure}

\begin{figure}[H]
\includegraphics{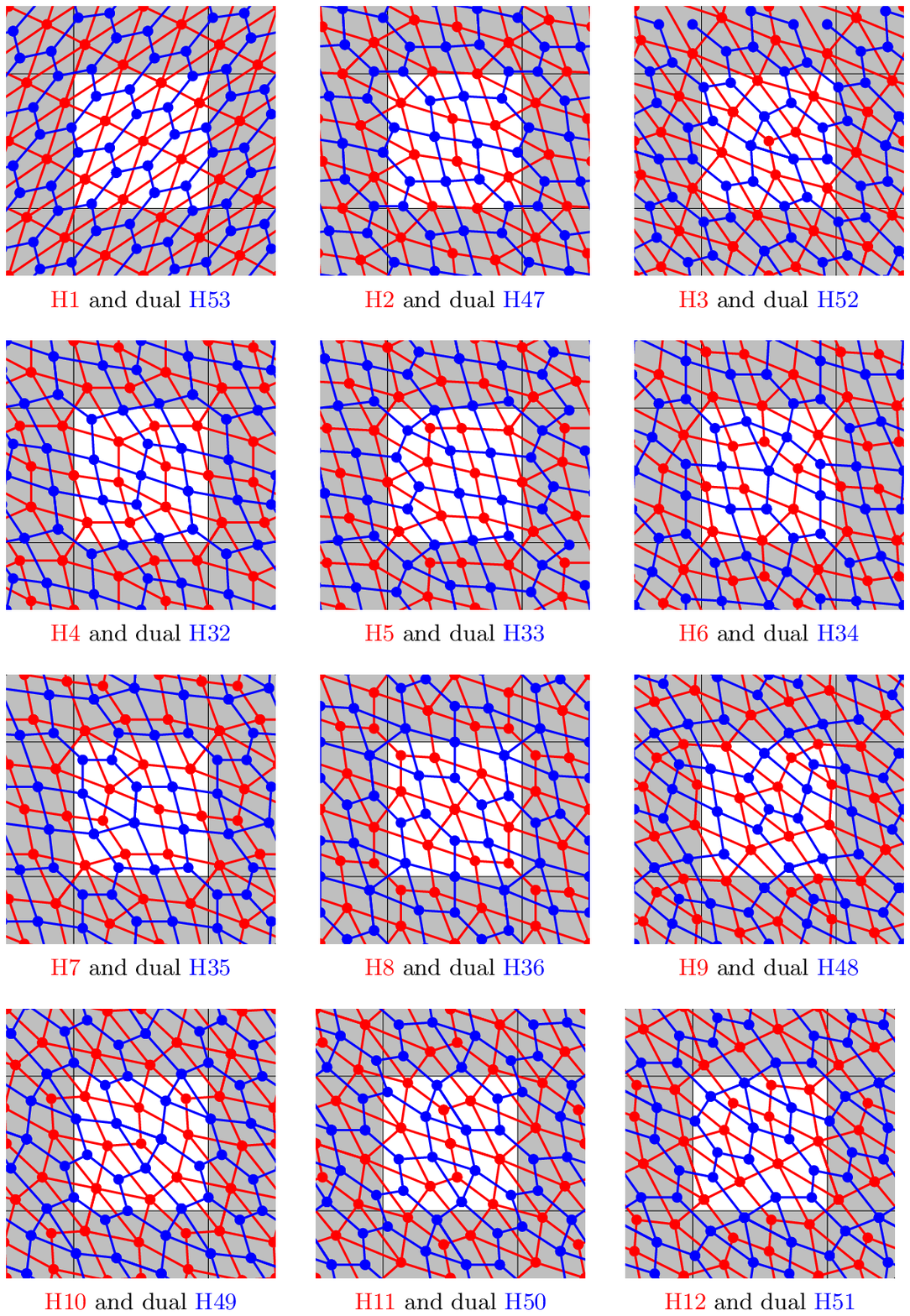}
\vspace*{3mm}
\caption{Diminimal Maps on the Torus with duals, 1 of 3}
\label{fig1}
\end{figure}

\begin{figure}[H]
\includegraphics{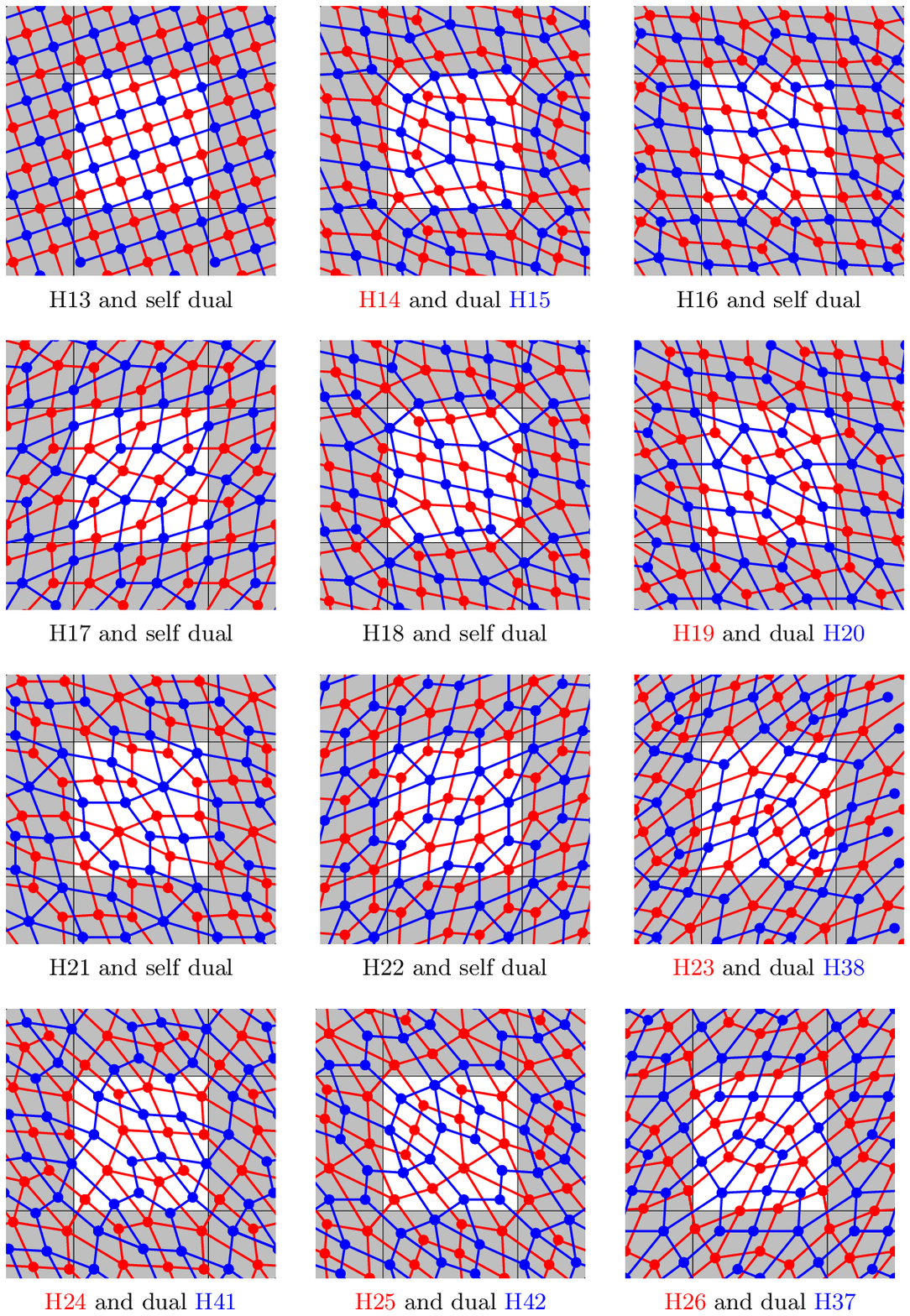}
\vspace*{3mm}
\caption{Diminimal Maps on the Torus with duals, 2 of 3}
\label{fig2}
\end{figure}

\begin{figure}[H]
\includegraphics{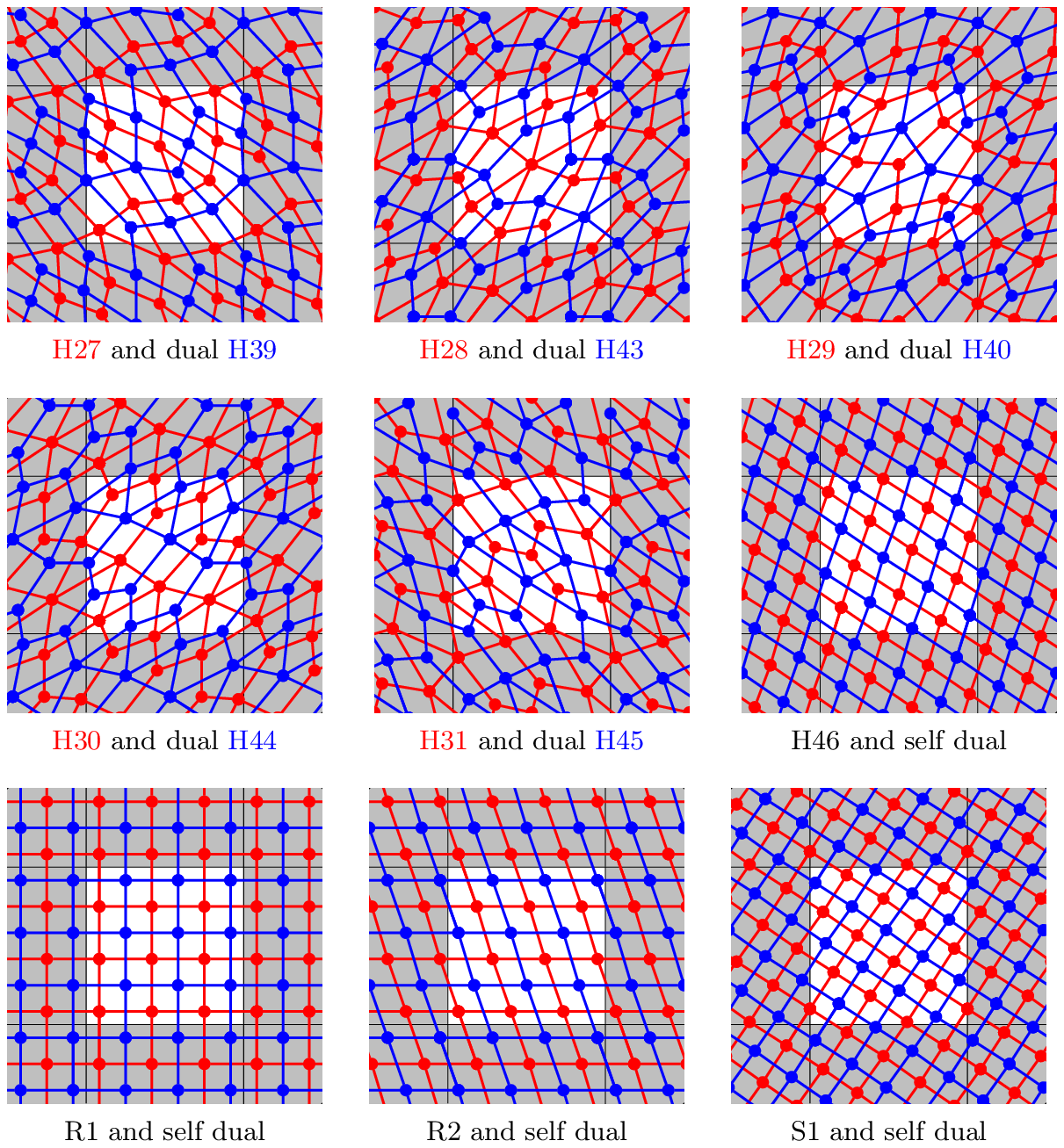}
\vspace*{3mm}
\caption{Diminimal Maps on the Torus with duals, 3 of 3}
\label{fig3}
\end{figure}


\bibliographystyle{amsplain}
\bibliography{diminimal}

\end{document}